\documentclass[11pt]{article}
\usepackage{amstext,amssymb,amsmath,amsbsy}

\usepackage{hyperref}
\usepackage{amscd}
\usepackage{amsfonts}
\usepackage{indentfirst} 
\usepackage{verbatim}
\usepackage{amsmath}
\usepackage{amsthm}
\usepackage{enumerate}
\usepackage{graphicx} 
\usepackage{color}
\usepackage[OT1]{fontenc}
\usepackage[latin1]{inputenc}
\usepackage[english]{babel}
\usepackage{amssymb,MnSymbol}
\usepackage{tikz}
\usepackage{dsfont}
\newtheorem{theorem}{Theorem}
\newtheorem{lemma}{Lemma}

\newtheorem{proposition}{Proposition}

\newtheorem{remark}{Remark}

\newcommand{\proofend}{\hfill $\Box$ }
\newcommand{\mint}{\strokedint}

\newcommand{\loc}{_{loc}}

\newcommand{\mN}{\mathbb{N}}

\newcommand{\mR}{\mathbb{R}}
\newcommand{\mS}{\mathbb{S}}
\newcommand{\eps}{\varepsilon}

\newcommand{\dsp}{\displaystyle}
\newcommand{\mc}{\mathrm{c}}

\newcommand{\supp}{\operatorname{supp}}

\numberwithin{equation}{section}

\title{The BBM formula revisited}

\author{Ha\"im Brezis\footnote{Rutgers University,
Department of Mathematics, Hill Center, Busch Campus,
110 Frelinghuysen Road, Piscataway, NJ 08854, USA, brezis@math.rutgers.edu} \footnote{Department of Mathematics,
Technion, Israel Institute of Technology,
32.000 Haifa, Israel} \footnote{Laboratoire Jacques-Louis Lions
UPMC,  4  place Jussieu, 75005 Paris,
France} \footnote{Research partially supported by NSF grant DMS-1207793 and by ITN "FIRST" of the European Commission, Grant  Number  PITN-GA-2009-238702.} \; and Hoai-Minh Nguyen\footnote{EPFL SB MATHAA CAMA, Station 8,  CH-1015 Lausanne,  Switzerland, hoai-minh.nguyen@epfl.ch}
}

\begin{document}

\maketitle

\vspace{-0.7cm}
\begin{center}
{\it To the memory of Ennio De Giorgi with emotion and admiration}
\end{center}

\begin{abstract}
In this paper, we revise the BBM formula due to J. Bourgain, H. Brezis, and P. Mironescu in \cite{BBM}. 
\end{abstract}


\noindent Scientific chapter:  12. Real Variable(s) Functions. 



\section{Introduction}
We first recall the BBM formula due to J. Bourgain, H. Brezis, and P. Mironescu \cite{BBM}, see also \cite{B1},  (with a refinement by J. Davila \cite{Davila}).
Let $d \ge 1$ be an integer. Throughout this paper, $(\rho_n)$ denotes a sequence of radial mollifiers in the sense that
\begin{equation}
\label{rho-1} \rho_n \in L^1_{\loc}(0, + \infty), \quad \rho_n \ge 0,
\end{equation}
\begin{equation}
\label{rho-2} \int_{0}^\infty \rho_n(r) r^{d-1} \, dr = 1 \quad \forall \, n,
\end{equation}
and
\begin{equation}
\label{rho-3} \lim_{n \to + \infty} \int_\delta^\infty \rho_n(r) r^{d-1} \, dr = 0  \quad \forall \, \delta > 0.
\end{equation}
Even though the next assumption is required only  for a few  results,   it is convenient to assume that
\begin{equation}
\label{rho-4} \rho_n(r) = 0 \mbox{ for all } r > 1, n \in \mN.
\end{equation}
Set, for $p \ge 1$,
\begin{equation}
I_{n, p} (u) = \int_{\mR^d} \int_{\mR^d}  \frac{|u(x) - u(y)|^p}{|x - y|^p} \rho_n(|x-y|) \, dx \, dy \le +\infty, \quad \forall \, u  \in L^1_{\loc}(\mR^d).
\end{equation}
For $u \in L^1_{\loc}(\mR^d)$, define,  for $p> 1$,
\begin{equation}\label{def-I}
I_{p}(u) = \left\{ \begin{array}{cl} \dsp \gamma_{d, p} \int_{\mR^d} |\nabla u|^p &  \mbox{ if } \nabla u \in L^{p}(\mR^d), \\[6pt]
+ \infty & \mbox{ otherwise},
\end{array}\right.
\end{equation}
and,  for $p=1$,
\begin{equation}\label{def-I}
I_{1}(u) = \left\{ \begin{array}{cl} \dsp \gamma_{d, 1} \int_{\mR^d} |\nabla u| &  \mbox{ if } \nabla u \mbox{ is a finite measure}, \\[6pt]
+ \infty & \mbox{ otherwise},
\end{array}\right.
\end{equation}
where, for any $e \in \mS^{d-1}$ and $p \ge 1$,
\begin{equation}\label{def-gamma}
\gamma_{d, p} = \int_{\mS^{d-1}} |\sigma \cdot e|^p \, d \sigma.
\end{equation}
In the case $p=1$, we have
\begin{equation}\label{def-gamma}
\gamma_{d, 1} = \int_{\mS^{d-1}} |\sigma \cdot e| \, d \sigma = \left\{ \begin{array}{cl} \dsp \frac{2}{d-1} |\mS^{d-2}|=2 |B^{d-1}| & \mbox{ if } d \ge 3, \\[6pt]
4 & \mbox{ if } d = 2, \\[6pt]
2 & \mbox{ if } d =1.
\end{array} \right.
\end{equation}
The BBM formula asserts that, for $p \ge 1$,
\begin{equation}\label{limit-In}
\lim_{n \to + \infty} I_{n, p} (u)  = I_p(u) \quad  \forall \, u \in L^1_{loc}(\mR^d).
\end{equation}
Applying \eqref{limit-In} with $p=1$, $u = \mathds{1}_E$ (the characteristic function of a measurable set $E$), and $\rho_n (r) = C_d n^{(d+1)/2} r e^{-n r^2}$, we obtain
$$
\lim_{n \to + \infty} n^{(d+1)/2} \int_{E^c} \int_E e^{-n |x-y|^2} \, dx \, dy = A_d \mathrm{Per}(E).
$$
By comparison the De Giorgi formula \cite{DeGiorgi1, DeGiorgi2} for the perimeter involves a derivative and asserts that
$$
\lim_{n \to + \infty} \int_{\mR^d} |\nabla W_n(x)| \, dx = B_d \mathrm{Per} (E),
$$
where
$$
W_n(x) = n^{d/2} \int_{E} e^{- n |x - y|^2} \, dy,
$$
and $A_d$, $B_d$, and $C_d$ are positive constants depending only on $d$.

Define, for $p \ge 1$, $n \in \mN$, and $u \in L^1_{\loc}(\mR^d)$,
\begin{equation}
D_{n, p} (u) (x): = \int_{\mR^d}  \frac{|u(x) - u(y)|^p}{|x - y|^p} \rho_n(|x-y|) \, dy \mbox{ for a.e. } x \in \mR^d.
\end{equation}
Note that, see \cite{BBM},
$$
\int_{\mR^d} D_{n, p}(u) (x) \, dx  \le C_{p, d} \int_{\mR^d} |\nabla u|^p(x) \, dx \mbox{ for } n \in \mN,
$$
and hence
\begin{equation}\label{Dn-finite}
D_{n, p}(x) < + \infty \mbox{ for a.e. } x \in \mR^d
\end{equation}
if $p>1$ and  $\nabla u \in L^p(\mR^d)$ or $p=1$ and $ \nabla u $ is a finite measure.
From the BBM formula, we have, for $p \ge 1$,
\begin{equation}
\lim_{n \to + \infty} \int_{\mR^d} D_{n, p}(u) (x) = I_p(u) \quad \mbox{ for } u \in L^1_{\loc}(\mR^d).
\end{equation}
On the other hand, an easy computation (see \cite[formula (6)]{BBM}) gives, for $p \ge 1$, $u \in C^1_{\mc}(\mR^d)$, and $x \in \mR^d$,
$$
\lim_{n \to \infty} D_{n, p} (u)(x) =  \gamma_{d, p} |\nabla u|^p (x).
$$

In this paper, we investigate the mode convergence of $D_{n, p}(u)$ to $ \gamma_{d, p} |\nabla u|^p$ as $n \to + \infty$ for non smooth $u$.  Our main results are the following
\begin{theorem}\label{thm1} Let  $d \ge 1$, $p \ge 1$, and  $u \in W^{1, p}_{\loc}(\mR^d)$. Then
\begin{equation}\label{limit-ae}
\lim_{n \to + \infty} \int_{\mR^d} \frac{|u(x+h) - u(x) - \nabla u (x) \cdot h |^p}{|h|^p}\rho_n(|h|) \, dh   = 0 \mbox{ for  a.e. } x \in \mR^d.
\end{equation}
Consequently,
\begin{equation}\label{limit-ae-1}
\lim_{n \to + \infty} D_{n, p}(u)(x) = \gamma_{d, p} |\nabla u|^p(x) \mbox{ for a.e. } x \in \mR^d.
\end{equation}
\end{theorem}

\begin{remark}  \fontfamily{m} \selectfont When $\rho_n(r) = d \eps_n^{-d} \mathds{1}_{(0, \eps_n)}$ for a sequence of $(\eps_n) \to 0_+$, assertion \eqref{limit-ae} is part of the classical $L^{p}$-differentiability theory of Calder\'{o}n-Zygmund; the same comment applies to assertion \eqref{limit-ae-thm2} below. Theorem~\ref{thm1} is due to D. Spector \cite[Theorem 1.7]{S2} under  the additional assumption that $\rho_n$ is non-increasing for every $n$.  His argument is much more complicated than ours (in addition he relies on the $L^{p^*}$-differentiability of $W^{1, p}$ functions, see  e.g.,  \cite[Theorem 2 on page 262]{EGMeasure}).
\end{remark}

We now turn to the $L^1$-convergence of $D_{n, p}$.

\begin{proposition} \label{pro1-1}
Let  $d \ge 1$, $p \ge 1$, and  $u \in L^1_{\loc}(\mR^d)$ with $\nabla u \in L^p(\mR^d)$. Then
\begin{equation}\label{limit-L^p}
\lim_{n \to + \infty} \int_{\mR^d} \int_{\mR^d} \frac{|u(x+h) - u(x) - \nabla u (x) \cdot h |^p}{|h|^p}\rho_n(|h|) \, dh \, dx   = 0.
\end{equation}
Consequently,
\begin{equation}\label{limit-L^p-1}
\lim_{n \to + \infty} D_{n, p}(u)= \gamma_{d, p} |\nabla u|^p \mbox{ in } L^1(\mR^d).
\end{equation}
\end{proposition}

\begin{remark}  \fontfamily{m} \selectfont 
Assertion \eqref{limit-L^p-1} was proved in \cite{BBM}.  
\end{remark}

Theorem~\ref{thm1} (resp. Proposition~\ref{pro1-1}) is established in Section~\ref{sect-ae} (resp. Section~\ref{sect-L^1}) where we also present some variants, generalizations, and pathologies related to  these results. 

\medskip
The case $p=1$ and $u \in BV_{\loc}(\mR^d)$ is more delicate. In this case instead of Theorem~\ref{thm1}, we have
\begin{theorem}\label{thm2} Let  $d \ge 1$ and  $u \in BV_{\loc}(\mR^d)$. Then
\begin{equation}\label{limit-ae-thm2}
\lim_{n \to + \infty} \int_{\mR^d} \frac{|u(x+h) - u(x) - \nabla^{ac} u (x) \cdot h |}{|h|}\rho_n(|h|) \, dh   = 0 \mbox{ for  a.e. } x \in \mR^d.
\end{equation}
Consequently,
\begin{equation}\label{limit-ae-1-thm2}
\lim_{n \to + \infty} D_{n, 1}(u)(x) = \gamma_{d, 1} |\nabla^{ac} u|(x) \mbox{ for a.e. } x \in \mR^d.
\end{equation}
\end{theorem}

Here  and in what follows, for $u \in BV_{\loc}(\mR^d)$, we denote $\nabla^{ac} u$ and $\nabla^s u$ the absolutely continuous part and the singular part of $\nabla u$.

\begin{remark}  \fontfamily{m} \selectfont
 A version of Proposition~\ref{pro1-1}  for $u \in BV(\mR^d)$ has been established by A.~Ponce and D.~Spector \cite[Proposition 2.1]{PS}. Here is their result:   Let  $d \ge 1$,  and  $u \in BV(\mR^d)$. Then
\begin{equation*}
\lim_{n \to + \infty}  \int_{\mR^d} \frac{|u(x+h) - u(x) - \nabla^{ac} u (x) \cdot h| }{|h|}\rho_n(|h|) \, dh   = \gamma_{d, 1} |\nabla^s u| \mbox{ in the sense of measures}.
\end{equation*}
\end{remark}

Theorem~\ref{thm2} is established in Section~\ref{sect-BV}. In the last section, we present miscellaneous facts related to the above results.

\section{Convergence almost everywhere in the Sobolev case} \label{sect-ae}

We will use the following elementary lemma (see \cite[Lemma 1]{BrezisNguyen-Two}):

\begin{lemma} \label{lem-maximal-BN} Let $d \ge 1$,  $r>0$,  $x \in \mR^d$,  and $f \in L^1_{\loc}(\mR^d)$.  We have
\begin{equation}\label{maximal-part1}
\int_{\mS^{d-1}} \int_0^r |f(x + s \sigma)| \, d s \, d \sigma \le C_d r M(f)(x),
\end{equation}
for some positive constant $C_d$ depending only on $d$.
\end{lemma}

Here $M(f)$ denotes the maximal function of $f$.  We now give the

\medskip

\noindent{\bf Proof of Theorem~\ref{thm1}.}  We first present the proof for $u \in W^{1, p}(\mR^d)$.
We claim that, for all $u \in W^{1, p}(\mR^d)$,
\begin{equation}\label{bound}
D_{n, p}(u)(x) \le C M (|\nabla u|^p)(x) \mbox{ for  a.e. } x \in \mR^d.
\end{equation}
Here and in what follows, $C$ denotes a positive constant depending only on $d$.  We have, for a.e. $x \in \mR^d$, $\sigma \in \mS^{d-1}$, and $r > 0$,
$$
u(x + r \sigma ) - u(x) = \int_0^r \nabla u( x+ s \sigma ) \cdot \sigma \, ds.
$$
Using polar coordinates, H\"older's inequality,  and  Fubini's theorem,  we obtain, for a.e. $x \in \mR^d$,
\begin{align*}
\int_{\mR^d} \frac{|u(x+ h) - u(x)|^p}{|h|^p}\rho_n(|h|) \, dh \le &  \int_0^\infty \rho_n(r) r^{d-1} \frac{1}{r} \int_{\mS^{d-1}} \int_0^r |\nabla u(x+ s  \sigma) \cdot \sigma|^p \, d s \, d \sigma \, d r \\[6pt]
= & \int_0^\infty \rho_n(r) r^{d-1} \frac{1}{r} \int_{B(x,r)} |\nabla u(y)|^p |y|^{1-d} \, d y \, d r.
\end{align*}
Applying Lemma~\ref{lem-maximal-BN}, we obtain  \eqref{bound}.

The proof of \eqref{limit-ae} now goes  as follows. Set
$$
\Omega(u): = \Big\{x \in \mR^d;  \limsup_{n \to + \infty}  \int_{\mR^d} \frac{|u(x+h) - u(x) - \nabla u (x) \cdot h |^p}{|h|^p}\rho_n(|h|) \, dh  > 0 \Big\}.
$$
Note  that if $u \in C^1_{\mc}(\mR^d)$ then \eqref{limit-ae} holds for all $x \in \mR^d$. This implies
$$
|\Omega(v)| = 0 \mbox{ for all } v \in C^{1}_{\mc}(\mR^d).
$$
It follows that
\begin{equation}\label{state2}
\Omega(u) = \Omega(u - v) \mbox{ for all } v \in C^{1}_{\mc}(\mR^d).
\end{equation}
Recall that, see e.g.,  \cite[Theorem 1 on page 5]{SteinSingular}, for $f \in L^1(\mR^d)$, we have
\begin{equation}\label{theory-maximal}
\big|\{x \in \mR^d; M(f)(x) > \eps \} \big| \le \frac{C}{\eps} \int_{\mR^d} |f|.
\end{equation}
Using \eqref{bound} and \eqref{theory-maximal}, we obtain
\begin{multline}\label{state3}
\Big| \Big\{x \in \mR^d  \int_{\mR^d} \frac{|(u - v)(x+h) - (u-v)(x) - \nabla (u - v) (x) \cdot h |^p}{|h|^p}\rho_n(|h|) \, dh  > \eps \Big\}  \Big|\\[6pt]
 \le \frac{C}{\eps} \int_{\mR^d} |\nabla (u - v)(x)|^p \,dx \quad \mbox{ for all } \eps > 0.
\end{multline}
Combining  \eqref{state2} and \eqref{state3} yields \eqref{limit-ae}.   Assertion \eqref{limit-ae-1} follows from \eqref{limit-ae} by the triangle inequality after noting that, for every $V \in \mR^d$,
$$
\int_{\mR^d}\frac{ |V \cdot h|^p}{|h|^p} \rho_n(|h|) \, dh = \int_0^\infty \int_{\mS^{d-1}} |V \cdot \sigma |^p \rho_n(r) r^{d-1} \, d \sigma \, d r  = \gamma_{d, p} |V|^p.
$$

We now turn to the proof in the case $u \in W^{1, p}_{\loc}(\mR^d)$. Given $R> 1$, let $\varphi \in C^1_{\mc}(\mR^d)$ be such that $\varphi = 1$ in $B(0, 2 R)$. We have $\varphi u \in W^{1, p}(\mR^d)$. Applying the above result to $\varphi u$, we obtain
$$
\lim_{n \to + \infty} D_{n, p} (\varphi u)(x) = \gamma_{d, p} |\nabla (\varphi u)|^p(x) \mbox{ for a.e. } x \in B(0, R).
$$
Since $D_{n, p} (u)(x) = D_{n, p} (\varphi u)(x) $ for $x \in B_R$ by \eqref{rho-4} and $\varphi(x) u (x) = u(x)$ in $B_{R}$, it follows that
$$
\lim_{n \to + \infty}  D_{n, p} ( u)(x) = \gamma_{d, p} |\nabla ( u)|^p(x) \mbox{ for a.e. } x \in B(0, R).
$$
Since $R > 1$ is arbitrary, the conclusion follows. \proofend

\medskip
Here is a natural question related to Theorem~\ref{thm1}. Suppose for example that
$u \in W^{1,  1}(\mR^d)$ and  $u$ has compact support. Is it true that for every $1 < p < + \infty$,
$$
\lim_{n \to + \infty} D_{n, p}(u)(x) = \gamma_{d, p} |\nabla u|^p(x) \mbox{ for a.e. } x \in \mR^d?
$$
Surprisingly, the answer is delicate and some pathologies may occur as seen in our next result.

\begin{theorem} Let $d \ge 1$ and $u \in W^{1, 1}_{\loc}(\mR^d)$. We have
\begin{enumerate}
\item If $d=1$, then, for $p> 1$,
\begin{equation}\label{state1-Lemd=1}
\lim_{n \to + \infty} D_{n, p}(u)(x) = \gamma_{1, p} |u'|^p(x) \mbox{ for a.e. } x \in \mR.
\end{equation}

\item If $d \ge 2$, $p \le d/ (d-1)$,  and $\rho_n$ is non-increasing, then
\begin{equation}\label{state2-Lemd=1}
\lim_{n \to + \infty} D_{n, p}(u)(x) = \gamma_{d, p} |\nabla u|^p(x) \mbox{ for a.e. } x \in \mR^d.
\end{equation}

\item If $d \ge 2$ and  $p > 1$, then
\begin{equation}\label{state3-Lemd=1}
\liminf_{n \to + \infty} D_{n, p}(u)(x) \ge \gamma_{d, p} |\nabla u|^p(x) \mbox{ for a.e. } x \in \mR^d.
\end{equation}

Moreover, strict inequality in \eqref{state3-Lemd=1} can occur:

\item If $d \ge 2$, there exist $u \in W^{1, 1}(\mR^d)$ with compact support, a set $A \subset \mR^d$ of positive measure, and  a sequence of non-increasing functions $(\rho_n)$  such that, for every $n \in \mN$,
\begin{equation}\label{state4-Lemd=1}
D_{n, p} (u)(x) = + \infty \mbox{ for a.e. } x \in A, \quad \mbox{ for all } p> d/ (d-1).
\end{equation}
\end{enumerate}
\end{theorem}

Note that there is  no contradiction between \eqref{Dn-finite} and \eqref{state4-Lemd=1};  the $u$ which we construct here does not satisfy the condition   $\nabla u \in L^p(\mR^d)$.

\begin{remark} \label{rem-Spector}  \fontfamily{m} \selectfont  Statement \eqref{state2-Lemd=1} is due to D. Spector \cite[Theorem 1.7]{S2}.  In fact, he proves a more general result: if $u \in W^{1, q}(\mR^d)$ ($d \ge 2$) with $1 \le q < d$, $p\le  q^*= qd/ (d - q)$, and $\rho_n$ is non-increasing then \eqref{state2-Lemd=1} holds.
\end{remark}

\begin{remark}  \fontfamily{m} \selectfont  We do not know whether \eqref{state2-Lemd=1} holds without the additional assumption that $\rho_n$ is non-increasing.
\end{remark}

\noindent{\bf Proof.} As in the proof of Theorem~\ref{thm1}, one may assume that $u \in W^{1,1}(\mR^d)$.  We first prove \eqref{state1-Lemd=1}. Since, for a.e. $x \in \mR$ and $r > 0$,
$$
|u(x + r) - u(x)| \le \int_x^{x + r} |u'(s)| \, ds,
$$
we have
\begin{equation*}
D_{n,p}(u)^{1/p}(x) \le C M(u')(x).
\end{equation*}
Assertion \eqref{state1-Lemd=1} now follows as in the proof of Theorem~\ref{thm1} by noting that, for $u \in C^1_{\mc}(\mR)$,
$$
\lim_{n \to + \infty} D_{n, p}(u)(x) = \gamma_{1, p} |u'|^p(x) \mbox{ for } x \in \mR^d.
$$

We next turn to the proof of \eqref{state3-Lemd=1}. Using polar coordinates, we have, for a.e. $x \in \mR^d$,
\begin{align}\label{part1-state3}
D_{n, p}(u)(x) &  =  \int_{0}^\infty \int_{\mS^{d-1}} \left| \int_0^1 \nabla u(x + t r \sigma) \cdot \sigma \, d t \right|^p \rho_n(r) r^{d-1} \, d \sigma  \, dr \nonumber \\[6pt]
& \ge \int_{\mS^{d-1}} \left| \int_{0}^\infty\int_0^1 \nabla u(x + t r \sigma) \cdot \sigma \,    \rho_n(r) r^{d-1}  \, dt \, dr \right|^p \, d \sigma.
\end{align}
We claim that, for a.e. $\sigma \in \mS^{d-1}$ and for a.e. $x \in \mR^d$,
\begin{equation}\label{part2-state3}
\lim_{n \to + \infty} \int_{0}^\infty\int_0^1 \nabla u(x + t r \sigma) \cdot \sigma \,    \rho_n(r) r^{d-1}  \, dt \, dr  = \nabla u(x) \cdot \sigma.
\end{equation}
Assuming this and applying Fatou's lemma, we derive from \eqref{part1-state3} and \eqref{part2-state3} that, for a.e. $x \in \mR^d$,
\begin{equation*}
\liminf_{n \to + \infty}D_{n, p}(u)(x) \ge \gamma_{p, d} |\nabla u|^p(x);
\end{equation*}
which is \eqref{state3-Lemd=1}. To complete the proof of \eqref{state3-Lemd=1}, it remains to prove \eqref{part2-state3}. For $v \in W^{1, 1}(\mR^d)$,  $x \in \mR^d$,  and $\sigma \in \mS^{d-1}$, set
\begin{equation}
M(\nabla v, \sigma, x) = \sup_{r > 0} \mint_{0}^r |\nabla v (x + s \sigma ) \cdot \sigma | \, ds.
\end{equation}
Given $v \in W^{1, 1}(\mR^d)$ and $\sigma \in \mS^{d-1}$,  we claim that for all $\eps > 0$, there exists a positive constant $C$ independent of $v$, $\eps$, and $\sigma$ such that
\begin{equation}\label{claim1-lem1}
\Big| \big\{x \in \mR^d; M(\nabla v, \sigma, x) > \eps\;  \big\} \Big| \le \frac{C}{\eps} \int_{\mR^d} |\nabla v(y)| \, dy.
\end{equation}
Using Fubini's theorem, we derive from \eqref{claim1-lem1} that
\begin{equation}\label{claim2-lem1}
\Big| \big\{(x, \sigma) \in \mR^d \times \mS^{d-1}; M(\nabla v, \sigma, x) > \eps\;  \big\} \Big| \le \frac{C}{\eps} \int_{\mR^d} |\nabla v(y)| \, dy.
\end{equation}
Using \eqref{claim2-lem1}, one can now obtain assertion  \eqref{part2-state3} as in the proof of Theorem~\ref{thm1} by noting that for all $u \in C^1_{\mc}(\mR^d)$,
\begin{equation*}
\lim_{n \to + \infty} \int_{0}^\infty\int_0^1 \nabla u(x + t r \sigma) \cdot \sigma \,    \rho_n(r) r^{d-1}  \, dt \, dr  = \nabla u(x) \cdot \sigma \quad \mbox{ for all } x \in \mR^d.
\end{equation*}
We next establish \eqref{claim1-lem1}. For simplicity of notation, we assume that $\sigma = e_d: =  (0, \cdots, 0, 1)$. We have, by Fubini's theorem,
\begin{equation}\label{toto1-lem1}
\Big| \big\{x \in \mR^d; M(\nabla v, e_d, x) > \eps\;  \big\} \Big| = \int_{\mR^{d-1}} \int_{\mR}  \mathds{1}_{\big\{x \in \mR^d; M(\nabla v, e_d, x) > \eps\;  \big\}} \, dx_d \, dx'.
\end{equation}
It follows from the theory of maximal functions (see \eqref{theory-maximal}) that
\begin{equation}\label{toto2-lem1}
\int_{\mR^{d-1}} \int_{\mR}  \mathds{1}_{\big\{x \in \mR^d; M(\nabla v, e_d, x) > \eps\;  \big\}} \, dx_d \, dx' \le \frac{C}{\eps} \int_{\mR^{d-1}} \int_{\mR} |\partial_{x_d} v (x', x_d)| \, dx_d \, d x' .
\end{equation}
Combining \eqref{toto1-lem1} and \eqref{toto2-lem1} yields
\begin{equation*}
\Big| \big\{x \in \mR^d; M(\nabla v, e_d, x) > \eps\;  \big\} \Big|  \le \frac{C}{\eps} \int_{\mR^{d}} |\nabla v(x)| \, dx;
\end{equation*}
which is \eqref{claim1-lem1}. The proof of \eqref{state3-Lemd=1} is complete.

\medskip

We finally establish \eqref{state4-Lemd=1}.  Let $(\delta_n) $ be a positive sequence converging to 0 such that $\delta_n < 1/2$ for all $n$, and define
\begin{equation}\label{pro-patho}
\rho_n(t) = \delta_n t^{\delta_n - 1} \mathds{1}_{(0, 1)} (t).
\end{equation}
Set  $u(x) =  \varphi(x) |x|^{(1-d)} \ln^{-2} |x| $ for some $\varphi \in C^1_{\mc}(\mR^d)$ such that $\varphi (x) = 1$ for $|x| < 2$. It is clear  that $u \in W^{1, 1}(\mR^d)$ and for $x \in \mR^d$ with $1/ 4 < |x| < 1/ 2$,
\begin{equation*}
\mathop{\int}_{|y | < 1/8}|u(x) - u(y)|^p \, dy = + \infty
\end{equation*}
since $p> d/ (d-1)$ and $\rho_n(|y-x|) \ge \delta_n (1/8)^{\delta_n -1}$ for $|y| < 1/8$ and $1/4 < |x| < 1/2$.
It follows that, for $1/4 < |x| < 1/2$,
\begin{equation*}
D_{n, p}(u)(x) = + \infty \quad \forall \, n.
\end{equation*}
The proof is complete.
\proofend

\section{Convergence in norm} \label{sect-L^1}

We present two proofs of Proposition~\ref{pro1-1}.

\medskip
\noindent{\bf First proof of Proposition~\ref{pro1-1} via Theorem~\ref{thm1}.}  By Theorem~\ref{thm1}, we have
\begin{equation}\label{state1-cor1}
\lim_{n \to + \infty  }D_{n, p }(u)(x) = \gamma_{d, p}  |\nabla u(x)|^p  \mbox{ for a.e. } x \in \mR^d.
\end{equation}
On the other hand, by the BBM formula,
\begin{equation}\label{state2-cor1}
\lim_{n \to + \infty} \int_{\mR^d} D_{n, p}(u)(x) \, dx =  \gamma_{d, p} \int_{\mR^d} |\nabla u(x)|^p \, dx.
\end{equation}
Recall that (see e.g., \cite[page 113]{BrAnalyseEnglish}) if $f_n(x) \to f(x)$ for a.e. $x \in \mR^d$,  and
$\| f_n\|_{L^1(\mR^d)} \to \| f\|_{L^1(\mR^d)}$,  then $f_n \to f$ in $L^1(\mR^d)$.
We  deduce from \eqref{state1-cor1} and \eqref{state2-cor1} that
$$
D_{n, p} (u) \to \gamma_{d, p} |\nabla u|^p \mbox{ in  } L^1(\mR^d)  \mbox{ as } n \to + \infty.
$$
\proofend

\medskip
\noindent{\bf Direct proof of Proposition~\ref{pro1-1}.}  We have, see \cite{BBM},
\begin{equation*}
\int_{\mR^d} \int_{\mR^d}  \frac{|u(x + h) - u(x) - \nabla u(x) \cdot  h |^p}{|h|^p} \rho_n(|h|) \, dh \, dx \le C_{p, d} \int_{\mR^d} |\nabla u (x)|^p
 \end{equation*}
and, for $v \in C^{1}_{\mc}(\mR^d)$,
\begin{equation*}
\lim_{n \to + \infty } D_{n, p }(v)(x) = \gamma_{d, p}  |\nabla v(x)|^p
 \mbox{ in } L^1(\mR^d) \mbox{ as } n \to + \infty.
\end{equation*}
The conclusion now follows by a standard approximation argument.  \proofend


\section{Convergence almost everywhere in the BV case} \label{sect-BV}

Let $d \ge 1$,  $\mu$ be a  Radon measure defined on $\mR^d$,  and $0 < R \le + \infty$. Denote
$$
M_R(\mu)(x)= \sup_{0 < s \le  R} \frac{|\mu|(B(x, s))}{|B(x, s)|} \quad \mbox{ and } \quad M(\mu)(x) = M_{\infty}(\mu)(x).
$$
We begin this section with

\begin{lemma}\label{lem-maximal}  Let  $d \ge 1$,  $\mu$ be a positive Radon measure defined in $\mR^d$, and let
$(\chi_k)_{k \ge 1}$ be a sequence of mollifier such that $\supp \chi_k \subset B(0, 1/k)$ and $0 \le \chi_k \le C k^d$ for some positive constant  $C$ depending only on $d$.  Set
$\mu_k = \mu* \chi_k$.  We have, for $x \in \mR^d$ and for $r>0$,
\begin{equation}\label{lem-maximal-1}
 \frac{1}{r} \int_{B(x, r)} |y - x|^{1- d} \, d \mu (y) \le C M_r(\mu)(x)
\end{equation}
and,  for every $k$,
\begin{equation}\label{lem-maximal-2}
 \frac{1}{r} \int_{B(x, r)} |y - x|^{1- d} \, d \mu_k (y) \le C M(\mu)(x),
\end{equation}
for some positive constant $C$ depending only on $d$.
\end{lemma}

\noindent{\bf Proof.} Without loss of generality, one may assume that $x =0$. We have
\begin{align*}
 \frac{1}{r} \int_{B(0, r)} |y|^{1- d} \, d \mu (y) &=  \frac{1}{r} \sum_{m = 0}^\infty \int_{B(0, 2^{-m} r) \setminus B(0, 2^{-(m+1)} r) } |y|^{1- d} \, d \mu (y) \\[6pt]
& \le   \frac{C}{r} \sum_{m = 0}^\infty  2^{-m(1-d)} r^{1-d} \int_{B(0, 2^{-m} r) \setminus B(0, 2^{-(m+1)} r) }   \, d \mu (y) \\[6pt]
& \le   \frac{C}{r} \sum_{m = 0}^\infty 2^{-m} r M_r(\mu)(0) = C M_r(\mu)(0);
\end{align*}
which is \eqref{lem-maximal-1}.

We next prove \eqref{lem-maximal-2}. As above, we obtain
\begin{equation}\label{max-muk-1}
 \frac{1}{r} \int_{B(0, r)} |y|^{1- d} \, d \mu_k (y) \le  \frac{C}{r} \sum_{m = 0}^\infty  2^{-m(1-d)} r^{1-d} \int_{B(0, 2^{-m} r) \setminus B(0, 2^{-(m+1)} r) }   \, d \mu_k (y).
\end{equation}
We claim that
\begin{equation}\label{max-muk-2}
\int_{B(0, 2^{-m} r) \setminus B(0, 2^{-(m+1)} r) }   \, d \mu_k (y) \le C 2^{-m d} r^d M(\mu)(0).
\end{equation}
Combining \eqref{max-muk-1} and \eqref{max-muk-2} yields \eqref{lem-maximal-2}

It remains to prove \eqref{max-muk-1}.  We have
\begin{align}\label{max-muk-3}
\int_{B(0, 2^{-m} r) \setminus B(0, 2^{-(m+1)} r)}\, d \mu_k (y) \le &  \int_{B(0, 2^{-m} r) \setminus \overline{B(0, 2^{-(m+2)} r)} }\, d \mu_k (y) \nonumber \\[6pt]
 =  & \sup_{\varphi \in C_{\mc} \big(B(0, 2^{-m} r) \setminus \overline{B(0, 2^{-(m+2)} r) } \big); |\varphi| \le 1} \int_{\mR^d} \varphi  \, d \mu_k.
\end{align}
We have
\begin{equation}\label{max-muk-4-0}
\int_{\mR^d} \varphi  \, d \mu_k =  \int_{\mR^d }\int_{\mR^d} \varphi (z) \chi_k(z - y) \, dz   \, d \mu (y)
\end{equation}
If $2^{-m} r < 1/k$, we have,  for $\varphi \in C_{\mc} \big(B(0, 2^{-m} r) \setminus \overline{B(0, 2^{-(m+2)} r)} \big)$ with $|\varphi| \le 1$,
\begin{align}\label{max-muk-4}
 \int_{\mR^d }\int_{\mR^d} \varphi (z) \chi_k(z - y) \, dz   \, d \mu (y) & \le  \int_{|y| < 2/ k} \sup_{y}\int_{\mR^d} |\varphi (z) | \chi_k(z - y) \, dz   \, d \mu (y)  \nonumber \\[6pt]
& \le  C (2^{-m} r)^d k^d  \int_{|y| < 2/ k}   \, d \mu (y)  \le C 2^{-md} r^d M(\mu)(0).
\end{align}
Here we use the fact that $\supp \chi_k \subset B(0, 1/k)$ and $0 \le \chi_k \le C k^d$.
Similarly, if $1/k < 2^{-m} r$, we have,  for $\varphi \in C_{\mc} \big(B(0, 2^{-m} r) \setminus \overline{B(0, 2^{-(m+2)} r)} \big)$ with $|\varphi| \le 1$,
\begin{align}\label{max-muk-5}
 \int_{\mR^d }\int_{\mR^d} \varphi (z) \chi_k(z - y) \, dz   \, d \mu (y) \, dy & \le \int_{|y| < 2^{-m + 2} r} \sup_{y}\int_{\mR^d} |\varphi (z) | \chi_k(z - y) \, dz   \, d \mu (y) \nonumber  \\[6pt]
&\le  \int_{|y| < 2^{-m + 2} r}  \, d \mu (y)   \le C 2^{-md} r^d M(\mu)(0).
\end{align}
Combining \eqref{max-muk-3}, \eqref{max-muk-4-0}, \eqref{max-muk-4}, and \eqref{max-muk-5}, we obtain \eqref{max-muk-2}. The proof is complete. \proofend

\medskip
We recall that  (see,  e.g.,  \cite{EGMeasure})
\begin{equation}\label{state1-BV1}
\lim_{r \to 0} \frac{|\nabla^{s} u| \big(B(x, r) \big)}{|B(x, r)|} = 0 \mbox{ for a.e. } x \in \mR^d.
\end{equation}
As a consequence of  \eqref{state1-BV1}, one obtains
\begin{equation}\label{state2-BV1}
M(|\nabla^s u|)(x) < + \infty \mbox{ for a.e. } x \in \mR^d.
\end{equation}

We now present the

\medskip
\noindent{\bf Proof of Theorem~\ref{thm2}.}  As in the proof of Theorem~\ref{thm1}, one may assume that $u \in BV(\mR^d)$.
Let $(\chi_k)_{k \ge 1}$ be a sequence of  smooth mollifiers such that $\supp \chi_k \subset B(0, 1/ k)$ and $0 \le \chi_k \le Ck^d$. Here and in what follows, $C$ denotes a positive constant depending only on $d$. Set, for $k \in \mN_+$,
$$
u_k = u* \chi_k, \quad V^s_k = \nabla^s u * \chi_k, \quad \mbox{ and } \quad V^{ac}_k = \nabla^{ac} u* \chi_k.
$$
We have
\begin{multline}\label{state3-BV1}
 \int_{\mR^d} \frac{|u_k(x+h) - u_k(x) - V^{ac}_k (x) \cdot h |}{|h|}\rho_n(|h|) \, dh \\[6pt]
 = \int_0^\infty r^{d-1} \rho_n (r) \int_{\mS^{d-1}} \frac{|u_k(x+r \sigma ) - u_k(x) - r V^{ac}_k (x) \cdot  \sigma |}{r} \, d \sigma \, d r.
\end{multline}
Since
\begin{equation*}
u_k(x+r \sigma ) - u_k(x) - r V^{ac}_k (x) \cdot  \sigma  = \int_0^{r} \nabla u_k (x + s \sigma)  \cdot \sigma \, ds  - r V^{ac}_k (x) \cdot \sigma
\end{equation*}
and
\begin{equation*}
\nabla u_k (x) = V^{s}_k(x) + V^{ac}_k(x),
\end{equation*}
it follows from \eqref{state3-BV1} that
\begin{multline}\label{state4-BV1}
\int_{\mR^d} \frac{|u_k(x+h) - u_k(x) - V^{ac}_k (x) \cdot h |}{|h|}\rho_n(|h|) \, dh \\[6pt]
\le  \int_0^\infty r^{d-1} \rho_n (r) \frac{1}{r}  \, dr \int_{\mS^{d-1}}   \int_0^r |V_k^{s}(x + s \sigma)| \, ds \, d \sigma \\[6pt]+  \int_0^\infty r^{d-1} \rho_n (r) \frac{1}{r} \, dr  \int_{\mS^{d-1}} \int_0^r |V_k^{ac}(x + s \sigma) - V_k^{ac} (x)|  \, ds \, d \sigma.
\end{multline}
We claim that, for a.e. $x \in \mR^d$,
\begin{multline}\label{claim1-BV1}
\lim_{k \to + \infty} \int_{\mR^d} \frac{|u_k(x+h) - u_k(x) - V^{ac}_k (x) \cdot h |}{|h|}\rho_n(|h|) \, dh  \\[6pt]
=  \int_{\mR^d} \frac{|u(x+h) - u(x) - \nabla^{ac} u (x) \cdot h |}{|h|}\rho_n(|h|) \, dh,
\end{multline}
\begin{multline}\label{claim2-BV1}
\lim_{k \to + \infty}  \int_0^\infty r^{d-1} \rho_n (r) \frac{1}{r} \, dr \int_{\mS^{d-1}} \int_0^r |V_k^{s}(x + s \sigma)| \, ds \, d \sigma  \\[6pt]
=  \int_0^\infty r^{d-1} \rho_n (r)  \frac{1}{r}  \, dr \int_{B(x, r)} |\nabla^s u(y)| |y - x|^{1 - d} \, dy,
\end{multline}
and
\begin{multline}\label{claim3-BV1}
\lim_{k \to + \infty}  \int_0^\infty r^{d-1} \rho_n (r) \frac{1}{r} \, dr \int_{\mS^{d-1}} \int_0^r |V_k^{ac}(x + s \sigma) - V_k^{ac} (x)|  \, ds \, d \sigma  \\[6pt]
= \int_0^\infty r^{d-1} \rho_n (r) \frac{1}{r} \, dr \int_{\mS^{d-1}} \int_0^r |\nabla ^{ac} u (x + s \sigma) - \nabla^{ac} u (x)|  \, ds \, d \sigma.
\end{multline}
Assuming these claims, we continue the proof. Combining \eqref{state4-BV1}, \eqref{claim1-BV1}, \eqref{claim2-BV1}, and \eqref{claim3-BV1} yields, for a.e. $x \in \mR^d$,
\begin{multline}\label{state5-BV1}
\int_{\mR^d} \frac{|u(x+h) - u(x) - \nabla^{ac} u (x) \cdot h |}{|h|}\rho_n(|h|) \, dh \\[6pt]
\le \int_0^\infty r^{d-1} \rho_n (r)  \frac{1}{r} \, dr \int_{B(x, r)} |\nabla^s u(y)| |y - x|^{1 - d} \, dy  \\[6pt]+  \int_0^\infty r^{d-1} \rho_n (r) \frac{1}{r} \, dr \int_{\mS^{d-1}} \int_0^r |\nabla^{ac} u(x + s \sigma) - \nabla^{ac}  u(x)|  \, ds \, d \sigma .
\end{multline}
Hence it suffices to prove that,  for a.e. $x \in \mR^d$,
\begin{equation}\label{state6-BV1}
\lim_{n \to + \infty} \int_0^\infty r^{d-1} \rho_n (r) \frac{1}{r} \, dr  \int_{B(x, r)} |\nabla^s u(y)| |y -x |^{1-d} \, dy = 0
\end{equation}
and
\begin{equation}\label{state7-BV1}
\lim_{n \to + \infty} \int_0^\infty r^{d-1} \rho_n (r) \frac{1}{r} \, dr  \int_{\mS^{d-1}} \int_0^r |\nabla^{ac} u(x + s \sigma) - \nabla^{ac}  u(x)|  \, ds \, d \sigma =0.
\end{equation}
Note that assertion \eqref{state7-BV1} holds for every $ x \in \mR^d$ if  $u \in C^{1}_{\mc}(\mR^d)$ and, by Lemma~\ref{lem-maximal},
$$
\int_0^\infty r^{d-1} \rho_n (r) \frac{1}{r}  \, dr \int_{\mS^{d-1}} \int_0^r |\nabla^{ac} u(x + s \sigma) - \nabla^{ac}  u(x)|  \, ds \, d \sigma \le C M(|\nabla^{ac} u|)(x).
$$
As in the proof of Theorem~\ref{thm1}, we have, for a.e. $x \in \mR^d$,
\begin{equation*}
\lim_{n \to + \infty} \int_0^\infty r^{d-1} \rho_n (r) \frac{1}{r}  \, dr \int_{\mS^{d-1}}  \int_0^r |\nabla^{ac} u(x + s \sigma) - \nabla^{ac}  u(x)|  \, ds \, d \sigma  =0;
\end{equation*}
which is \eqref{state7-BV1}.

We next establish \eqref{state6-BV1}. By Lemma~\ref{lem-maximal}, we have
\begin{equation*}
\frac{1}{r}\int_{B(x, r)} |\nabla^s u(y)| |y - x|^{1-d} \, dy \le C M_r(|\nabla^s u|) (x).
\end{equation*}
It follows from \eqref{state1-BV1} that
$$
\lim_{n \to + \infty} \int_0^\infty r^{d-1} \rho_n (r) \frac{1}{r} \, dr  \int_{B(x, r)} |\nabla^s u(y)| |y- x|^{1-d} \, dy = 0 \mbox{ for a.e. } x \in \mR^d;
$$
which is \eqref{state6-BV1}.

It remains to prove claims \eqref{claim1-BV1}, \eqref{claim2-BV1}, and \eqref{claim3-BV1}. We begin with claim~\eqref{claim1-BV1}.
We have
\begin{multline*}
\int_{\mR^d} \frac{|u_k(x+h) - u_k(x) - V^{ac}_k (x) \cdot h |}{|h|}\rho_n(|h|) \, dh \\[6pt]
= \int_{0}^\infty \rho_n(r) r^{d-1} \frac{1}{r} \, dr  \int_{\mS^{d-1}} |u_k(x+ r \sigma) - u_k(x) - r V^{ac}_k (x) \cdot \sigma | \, d \sigma.
\end{multline*}
Using  Lemma~\ref{lem-maximal}, we derive from \eqref{state4-BV1} that
\begin{equation*}
\frac{1}{r} \int_{\mS^{d-1}} |u_k(x+ r \sigma) - u_k(x) - r V^{ac}_k (x) \cdot \sigma | \, d \sigma \le C M(|\nabla u|)(x).
\end{equation*}
Since for a.e. $x \in \mR^d$,
\begin{multline*}
 \lim_{k \to + \infty}  \frac{1}{r}\int_{\mS^{d-1}} |u_k(x+ r \sigma) - u_k(x) - r V^{ac}_k (x) \cdot \sigma | \, d \sigma  \\[6pt]
 =  \frac{1}{r}\int_{\mS^{d-1}} |u(x+ r \sigma) - u(x) - r \nabla^{ac} u (x) \cdot \sigma | \, d \sigma  \mbox{ for  a.e. }  r > 0,
\end{multline*}
it follows from the  dominated convergence theorem that, for a.e. $x \in \mR^d$,
\begin{multline*}
\lim_{k \to + \infty} \int_{\mR^d} \frac{|u_k(x+h) - u_k(x) - V^{ac}_k (x) \cdot h |}{|h|}\rho_n(|h|) \, dh \\[6pt]
= \int_{\mR^d} \frac{|u(x+h) - u(x) - \nabla^{ac} u(x) \cdot h |}{|h|}\rho_n(|h|) \, dh;
\end{multline*}
which is \eqref{claim1-BV1}.

The proof of \eqref{claim3-BV1} follows similarly. We finally establish \eqref{claim2-BV1}. Fix $\tau > 0$ (arbitrary). We have
\begin{align}\label{claim2-BV1-part1}
\int_0^\infty r^{d-1} \rho_n (r) \frac{1}{r} \, dr \int_{\mS^{d-1}} \int_0^r &  |V_k^{s}(x + s \sigma)| \, ds \, d \sigma  \nonumber \\[6pt]
= &  \int_\tau^\infty r^{d-1} \rho_n (r) \frac{1}{r}  \, dr \int_{B(x, r) \setminus B(x, \tau) }   |V_k^{s}(y)| |y - x|^{1-d} \, dy \nonumber \\[6pt]
&+ \int_\tau^\infty r^{d-1} \rho_n (r) \frac{1}{r}\, dr  \int_{B(x, \tau)}   |V_k^{s}(y)| |y - x|^{1-d} \, dy  \nonumber \\[6pt]
& + \int_0^\tau r^{d-1} \rho_n (r) \frac{1}{r} \, dr \int_{B(x, r) }   |V_k^{s}(y)| |y - x|^{1-d} \, dy.
\end{align}
We have, for a.e. $r > 0$,
\begin{equation*}
\lim_{ k \to + \infty } \frac{1}{r} \int_{B(x, r) \setminus B(x, \tau) }   |V_k^{s}(y)| |y - x|^{1-d} \, dy  =  \frac{1}{r} \int_{B(x, r) \setminus B(x, \tau) }   |\nabla^{s} u (y)| |y- x|^{1-d} \, dy
\end{equation*}
and, by Lemma~\ref{lem-maximal},
\begin{equation*}
\frac{1}{r} \int_{B(x, r) \setminus B(x, \tau) }   |V_k^{s}(y)| |y - x |^{1-d} \, dy  \le  C M(|\nabla u|)(x).
\end{equation*}
It follows from the dominated convergence theorem that
\begin{multline}\label{claim2-BV1-part2}
\lim_{k \to + \infty } \int_\tau^\infty r^{d-1} \rho_n (r)   \frac{1}{r} \, dr \int_{B(x, r) \setminus  B(x, \tau) }   |V_k^{s}(y)| |y - x|^{1-d} \, dy  \\[6pt]
=  \int_\tau^\infty r^{d-1} \rho_n (r)   \frac{1}{r}  \, dr \int_{B(x, r) \setminus B(x, \tau) }   |\nabla^{s} u (y)| |y - x|^{1-d} \, dy.
\end{multline}
On the other hand, by Lemma~\ref{lem-maximal},
\begin{equation}\label{claim2-BV1-part3}
\int_\tau^\infty r^{d-1} \rho_n (r) \frac{1}{r} \, dr \int_{B(x, \tau)}   |V^{s}_k u(y)| |y - x|^{1-d} \, dy  \le C M(|\nabla u|)(x) \int_\tau^\infty r^{d-1} \rho_n (r) \tau/ r \, dr
\end{equation}
and
\begin{equation}\label{claim2-BV1-part4}
\int_0^\tau r^{d-1} \rho_n (r) \frac{1}{r}  \, dr \int_{B(x, r)}   |V_k^{s}(y)| |y  - x |^{1-d} \, dy  \le C M(|\nabla u|)(x) \int_0^\tau r^{d-1} \rho_n (r) \, dr.
\end{equation}
Since
$$
\lim_{\tau \to 0 } \Big( \int_\tau^\infty r^{d-1} \rho_n (r) \tau/ r \, dr  + \int_0^\tau r^{d-1} \rho_n (r) \, dr  \Big)= 0,
$$
we obtain \eqref{claim2-BV1} from \eqref{claim2-BV1-part1}, \eqref{claim2-BV1-part2}, \eqref{claim2-BV1-part3}, and \eqref{claim2-BV1-part4}. The proof is complete.  \proofend


\section{Miscellaneous results}

\subsection{On a characterization of $W^{1,1}(\mR^d)$}

The following result deals with a ``converse'' of Proposition~\ref{pro1-1}. It is due to  D. Spector in  \cite[Theorem 1.3]{S1}  and \cite[Theorem 1.4]{S2} in the case $\rho_n(r) = d \eps_n^{-d} \mathds{1}_{(0, \eps_n)}$ for a sequence of $(\eps_n) \to 0_+$ and  to A. Ponce and D. Spector \cite[Remark 5]{PS} for a general sequence $(\rho_n)$. The proof we present here is more direct.

\begin{proposition} Let  $d \ge 1$ and  $u \in L^1(\mR^d)$. Then $u \in W^{1, 1}(\mR^d)$  if and only if there exists $U \in [L^1(\mR^d)]^d$ such that
\begin{equation}\label{limit}
\lim_{n \to + \infty} \int_{\mR^d} \int_{\mR^d} \frac{|u(x+h) - u(x) - U(x) \cdot h|}{|h|}\rho_n(|h|) \, dh \, dx = 0.
\end{equation}
\end{proposition}

\noindent{\bf Proof.}  We already know that \eqref{limit} holds for $u \in W^{1, 1}(\mR^d)$  with $\nabla u = U$ by Proposition~\ref{pro1-1}. It remains to prove that if  \eqref{limit} holds, then  $u \in W^{1, 1}(\mR^d)$.  Let $(\chi_k)$ be a sequence of standard  mollifiers.  Define
\begin{equation*}
u_k = u* \chi_k \quad \mbox{ and } \quad  U_k = U* \chi_k.
\end{equation*}
We have
\begin{multline*}
\int_{\mR^d} \int_{\mR^d} \frac{|u_k(x+h) - u_k(x) - U_k(x) \cdot h| }{|h|}\rho_n(|h|) \, dh  \, dx \\[6pt]
= \int_{\mR^d} \int_{\mR^d} \Big| \int_{\mR^d} u(x + h - y) \chi_k (y) \, dy - \int_{\mR^d} u(x  - y) \chi_k (y) \, dy - \int_{\mR^d} U(x - y) \cdot h \chi_k (y) \, dy \Big| |h|^{-1}\rho_n(|h|) \, dh  \, dx .
\end{multline*}
This implies
\begin{multline*}
\int_{\mR^d} \int_{\mR^d} \frac{|u_k(x+h) - u_k(x) - U_k(x) \cdot h| }{|h|}\rho_n(|h|) \, dh  \, dx \\[6pt]
\le \int_{\mR^d} \int_{\mR^d}  \int_{\mR^d} \frac{\big| u(x + h - y) - u(x  - y) -U(x - y) \cdot h \big|}{|h|} \chi_k (y) \, dy \rho_n(|h|) \, dh  \, dx.
\end{multline*}
A change of variables gives
\begin{multline*}
\int_{\mR^d} \int_{\mR^d} \frac{|u_k(x+h) - u_k(x) - U_k(x) \cdot h| }{|h|}\rho_n(|h|) \, dh  \, dx \\[6pt]
\le\int_{\mR^d} \int_{\mR^d} \frac{|u (x+h) - u(x) - U(x) \cdot h| }{|h|}\rho_n(|h|) \, dh  \, dx .
\end{multline*}
We derive from \eqref{limit} that, for $k > 0$,
 \begin{equation*}
\lim_{n \to + \infty }\int_{\mR^d} \int_{\mR^d} \frac{|u_k(x+h) - u_k(x) - U_k(x) \cdot h| }{|h|}\rho_n(|h|) \, dh  \, dx = 0.
\end{equation*}
Since $u_k$ is smooth, we obtain
\begin{equation*}
U_k = \nabla  u_k.
\end{equation*}
As $k \to + \infty$, $u_k \to u$ and $U_k \to U$ in $L^1(\mR^d)$, so that $u \in W^{1, 1}(\mR^d)$ and $ \nabla u = U$. \proofend

\subsection{The bounded domain case}

Most of the above results hold when $\mR^d$ is replaced by a smooth bounded domain $\Omega$ of $\mR^d$. Define, for $p \ge 1$, $n \in \mN$, and $u \in L^1_{\loc}(\Omega)$,
\begin{equation}
D_{n, p}^\Omega (u) (x): = \int_{\Omega}  \frac{|u(x) - u(y)|^p}{|x - y|^p} \rho_n(|x-y|) \, dy \mbox{ for a.e. } x \in \Omega.
\end{equation}

Here is a typical result:
\begin{theorem}
Let  $d \ge 1$, $p \ge 1$ and  $u \in W^{1, p}(\Omega)$. Then
\begin{equation}\label{limit-ae-1-O}
\lim_{n \to + \infty} D_{n, p}^\Omega (u)(x) = \gamma_{d, p} |\nabla u|^p(x) \mbox{ for a.e. } x \in \Omega.
\end{equation}
\end{theorem}

\noindent{\bf Proof.} Let $\tilde u$ be an extension of $u$ to $\mR^d$ such that $\tilde u \in W^{1, p}(\mR^d)$. Let $\omega \subset \subset \Omega$.  We have, for $x \in \omega$,
\begin{equation}\label{p1-1}
 D_{n, p}^\Omega (u)(x) = D_{n, p}(\tilde u)(x) -  \int_{\mR^d \setminus \Omega} \frac{|\tilde u(x) - \tilde u(y)|}{|x-y|} \rho_n(|x-y|) \, dy.
\end{equation}
Applying  Theorem~\ref{thm1} to $\tilde u$, we have for a.e. $x \in \omega$,
\begin{equation}\label{p1-2}
\lim_{n \to + \infty} D_{n, p}(\tilde u)(x) = \gamma_{d, p} |\nabla \tilde u|^p(x) =  \gamma_{d, p} |\nabla  u|^p(x).
\end{equation}
Since $\omega$ is arbitrary, it suffices to prove that for a.e. $x \in \omega$,
\begin{equation}\label{p1-2}
\lim_{n \to + \infty} \int_{\mR^d \setminus \Omega} \frac{|\tilde u(x) - \tilde u(y)|}{|x-y|} \rho_n(|x-y|) \, dy = 0.
\end{equation}
Let $\varphi \in C^{1}(\mR^d)$ be such that $\varphi = 1$ in $\mR^d \setminus \Omega$ and $\varphi = 0$ in $\omega$. Applying Theorem~\ref{thm1} to $\varphi \tilde u$, we obtain, for a.e. $x \in \omega$,
\begin{equation}\label{p1-3}
\lim_{n \to + \infty} \int_{\mR^d \setminus \Omega} \frac{|\tilde u(y)|}{|x-y|} \rho_n(|x-y|) \, dy = 0.
\end{equation}
On the other hand, for a.e. $x \in \omega$,
\begin{equation}\label{p1-4}
\lim_{n \to + \infty} \int_{\mR^d \setminus \Omega} \frac{|\tilde u(x)|}{|x-y|} \rho_n(|x-y|) \, dy =|u(x)| \lim_{n \to + \infty} \int_{\mR^d \setminus \Omega} \frac{1}{|x-y|} \rho_n(|x-y|) \, dy  =  0
\end{equation}
Assertion \eqref{p1-2} now follows from \eqref{p1-3} and \eqref{p1-4}.
\proofend

\providecommand{\bysame}{\leavevmode\hbox to3em{\hrulefill}\thinspace}
\providecommand{\MR}{\relax\ifhmode\unskip\space\fi MR }
\providecommand{\MRhref}[2]{%
  \href{http://www.ams.org/mathscinet-getitem?mr=#1}{#2}
}
\providecommand{\href}[2]{#2}

\end{document}